\documentclass{amsart}
\usepackage{amssymb,latexsym}
\begin{document}

\title[Automorphism Groups of Groups of Order $8p^2$]{The Automorphism Groups \\ of the Groups of Order $8p^2$}

\author{Walter Becker}
\address{266 Brian Drive, Warwick, Rhode Island 02886}
\author{Elaine W. Becker}
\address{American Mathematical Society, Providence, Rhode Island}
\maketitle
\begin{abstract}

		The automorphism groups for the groups of orders
		$8p$ and $8p^2$ are given. The calculations were done
		using the programming language CAYLEY. Explicit
		presentations for both the groups of these orders
		and their automorphism groups are given.
\end{abstract}
\section{Introduction}

	In the early development of group theory much effort 
was devoted to the determination of the number of groups of a
certain type, e.g., the number of groups of degree 8 say, or
the number of groups of a certain order, e.g., 72, $16p$ ($p$ being
any odd prime), etc. This is not as active an endeavor as it
was about 100 years ago, but some work is still being done in
this area. (See references in Appendix 0 of \cite{1} and the book of
M. W. Short \cite{2}.) Some of these authors obtained explicit
presentations for all of the groups of orders such as $8p$, $16p$,
$16p^2$,\ldots\ for $p$ being an odd prime. There has been very little
work devoted to the explicit determination of the automorphism
groups of these groups. It seemed to the authors that if explicit
presentations for these groups could be found, then one might be 
able to do the same for their automorphism groups. This has in
fact turned out to be the case. \\

	Knowledge of the automorphism groups of groups is useful in
a variety of contexts; e.g., in the general problem of group 
extensions, Aut($G$) is one more invariant for the group $G$ which might
serve to distinguish one group from a second nonisomorphic one. 
There are other questions that could also be asked, e.g, what  
relation if any exists between the structure of the group $G$ and that of its
automorphism group, Aut($G$). \\

	Desktop computers are now being used in many different
areas where the tried and true methods of the field are paper
and pencil types of proof. Modern algebra, and especially group 
theory is one such field. The use of computers to solve problems
in group theory has come about because of the work of individuals
such as Drs. J. Cannon and J. Neubueser and their computing systems
(CAYLEY/MAGMA and GAP respectively) \cite{3}. The development of 
computing algorithms per say, as useful as these computing procedures
may be by themselves does not contribute to the solution of typical 
problems in abstract algebra. The intelligent use of the computer for 
providing a large number of examples and counterexamples can lead to the
formulation of conjectures (would-be theorems if you will) that point
the way to further progress in many areas. This work attempts to show 
how one can use computers to formulate theorems that deal with the 
automorphism groups of finite groups. Additional questions involving 
automorphism groups that may be investigated by means of computers 
can be found in the article by D. MacHale \cite{4}.\\

	The calculations reported upon below were done using the
programming language CAYLEY. For the groups of order $8p^2$, the
calculations were explicitly done for the primes $p=3$, 5 and 7, as 
well as for some higher orders, e.g., $p=17$, when the systematic 
behavior for these higher orders is not clear from the lower-order 
calculational results. The tables list the calculational results as 
a function of $p$, for any odd prime.  In some cases, for example 
cyclic groups, dihedral and dicyclic groups, the results 
were proven a long time ago using alternate, and 
more traditional, group-theoretic methods and can be found in many 
group theory textbooks. In other cases the general expressions 
given are deductions (or conjectures if you will) as to what the 
general behavior of these automorphism groups will be based upon the 
explicit calculations done for the cases of the smaller primes. Each 
entry in the tables below can be looked upon as a theorem to be 
proven in the more traditional theorem-proof approach to group 
theory problems.\\

	The calculational results to be presented below almost 
certainly can be determined by other non-computer methods, but 
explicit calculational methods, where they exist, are difficult to
find; e.g., the standard reference on the automorphism groups for 
abelian groups is a 1928 article in German by K. Shoda, and the
standard work on automorphism group towers is the 1939 article by
H. Wielandt \cite{6}, also in German. With the large number of 
group theory textbooks on the market it is surprising that most of 
them contain only a brief discussion of automorphism groups and 
even less on methods for their calculation \cite{5}.\\

	We hope that the general patterns found in our 
calculations of these automorphism groups will encourage others to 
develop other, more standard group-theoretic methods for the 
calculation of the automorphism groups of finite groups. The 
structural patterns found in the order $8p$  and $8p^2$ groups' 
automorphism groups are present in almost all of the other cases 
that the author has investigated. The results of these other 
investigations for orders $16p$, $16p^2$, $32p$, $8p^3$,\ldots\  will be 
published elsewhere. \\

	A major reason for this publication is to encourage others 
to consider the study of automorphism groups of finite groups as  
a field ripe with possibilities for producing very general theorems 
--- a structural theory of automorphism groups, if you will. The 
relationship between the group and its automorphism group is 
also a good area to look at, since in many cases the automorphism
groups [$Aut(C_p @ T$)] take the form $Hol(C_p) \times $ (invariant factor). The 
\textquotedblleft invariant factor\textquotedblright\ should be determinable from the group $T$ alone. \\

	In what follows we will need some information on the groups
of orders 8 and $8p$  as well as their automorphism groups. This
material will be summarized below for the reader's convenience.
Much of this material is available in standard books on group 
theory, e.g., Carmichael \cite{7} and Burnside \cite{8}, among others. Some 
of the material on the automorphism groups for the groups of
order $8p$  may be found in the books by Coxeter and Moser \cite{9},
and Wood and Thomas \cite{10}. We have not been able to obtain a 
complete \textquotedblleft closed form\textquotedblright\ expression for certain of the automorphism 
groups for the orders $8p^2$, namely those for the groups in which 
the action of the 2-group on the $p$-group is by a $Q_2$ image in which
$p\equiv 7$ mod(8). The problem is determining a representation of the
group $<2,3,4>$ of order 48 in the groups GL(2,$p$) for $p\equiv 7$ mod(8)
(see appendix 2). Section 3 contains a brief discussion of the 
automorphism groups of the groups of order $8p^2$ and some possible
generalizations. Section 4 contains additional material on the 
$C_8$ extensions. Certain questions on determining when two groups
are isomorphic or distinct may provide for interesting classroom
exercises in elementary group theory. One can certainly use some
fancy programming system (e.g., GAP) to do this, but it might be
useful to see why/how these various groups of order $8*17^2$ differ.  \\

\section{Background material on groups of order 8 and $8p$\\  and their
automorphism groups}
		
The five groups of order 8 and their automorphism groups are:\\

\begin{tabular}{|c|c|c|c|c|c|}\hline
     group $g$&$    C_8   $&$  C_4 \times  C_2 $&$  C_2 \times  C_2 \times  C_2 $&$   D_4 $&$  Q_2   $\\ \hline
    aut($g$) &$  C_2 \times  C_2 $&$    D_4  $& GL(3, 2) &$  D_4 $&$   S_4 $ \\ \hline
\end{tabular}
\linebreak
	
The group GL(3,2) is the simple group of order 168, and $S_4$ (the symmetric group of order 24) is a complete group.\\

	A complete group is a group all of whose automorphisms are inner. An alternate definition of a complete group is a group with a trivial center and whose automorphism group is isomorphic
to the group itself.\\

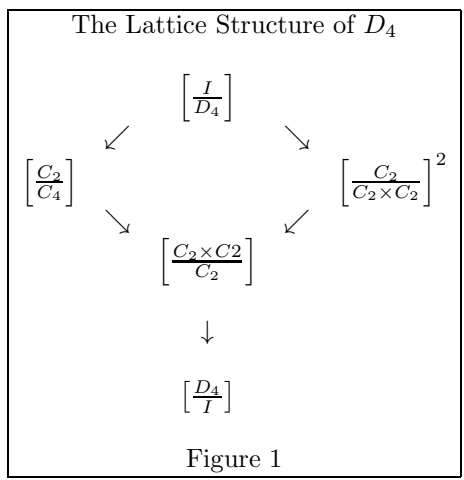
\begin{tabular}{|c c c c c |}\hline 
\multicolumn{5}{|c|}{The Lattice Structure of $D_4$}\\
&&&&\\  
& & $\left[\frac{I}{D_4}\right]$ & &  \\ 
& $\swarrow $& & $\searrow $&  \\ 
$\left[\frac{C_2}{C_4}\right]$ & & &  & $\left[\frac{C_2}{C_2\times C_2}\right]^2$ \\
& $ \searrow$&&$\swarrow$& \\ 
&&$\left[\frac{C_2 \times C2}{C_2}\right]$ & &  \\ 
&&&&\\
&& $\downarrow$ && \\
&&&&\\ 
& & $\left[\frac{D_4}{I}\right]$&& \\
&&&&\\ 
\multicolumn{5}{|c|}{Figure 1 }\\ \hline 
\end{tabular}\\
\linebreak
	To obtain the groups of order $8p$  where $p$ is an odd prime we need to know the normal subgroup structure of the groups of order 8. This information is readily available in many textbooks, e.g., the Hall-Senior tables \cite{11} and in \cite{10}. A typical example is the one for $D_4$ shown in Figure 1. From this figure we see that we have two different quotient groups of order 2, one with a kernel $C_4$ and the second one with the kernel $C_2 \times  C_2$. For the groups of order 8 there are seven such distinct normal subgroups, and each such $C_2$ quotient gives rise to a distinct group of order $8p$. All of these groups will appear in every order for which there exist groups of order $8p$. There are two groups of order 8 that have a quotient group $C_4$; these two groups will give rise to an additional pair of groups of order $8p$ for those primes that are equal to 1 mod(4). There is one additional group, coming from $C_8$, that will give rise to an additional group when $p\equiv 1$ mod(8). There are two groups of order 8 that have an automorphism of order 3, namely $C_2 \times  C_2 \times  C_2$ and $Q_2$. For the order 
24 these two groups will give rise to an additional pair of groups (namely $A_4 \times  C_2$ and SL(2,3)). In order 24 then we have the following groups:\\

\begin{tabular}{|c|} \hline
  5 direct products of the form $C_3 \times $ group of order 8  \\ \hline
  7 $C_2$-type extensions \\ \hline                                  
 2 groups having a normal sylow 2-subgroup: \\
                                              \\
   $ (A_4 \times  C_2$ and SL(2,3))         \\ \hline
  1 group without a normal sylow subgroup, $S_4$  \\ \hline   
\end{tabular}
\linebreak	

For order 40 ($p=5$ so $p$ is of the form 1 mod(4)), we have the five direct products plus the 7 $C_2$-type extensions plus the two coming from the $C_4$ action on $C_5$. There are no automorphisms of order 5 in the automorphism groups of order 8; so there are no groups of order 40 with a normal sylow 2-subgroup. All of the groups of order 40 have a normal sylow 5-subgroup. For $p = 7$ we have one extra group, namely ($C_2 \times  C_2 \times  C_2$) @ $C_7$, the Frobenius group of order 56.\\

	Presentations for the groups of order 8 are listed in Table 1. The non-direct-product groups of order $8p$  along with their automorphism group factors are also listed in Table 1. A simple illustration will show the reader the meaning of the entries in this table. Consider the extensions coming from the group $D_4$. From Table 1, $D_4$ has the presentation $a^4=b^2=a^b*a=1$. The group $D_4$ has two distinct normal subgroups of order 4 ($C_2 \times  C_2$ and $C_4$), giving rise to quotient groups of order 2. The two possible actions of $D_4$ on $C_p$ are:
\begin{equation}
\text{a)} \quad c^p=a^4=b^2=a^b*a=c^a*c=(b,c)=1\ \text{[case a in Table 1]}
\end{equation}
for the case of the $C_2 \times  C_2$ kernel and
\begin{equation}
\text{b)} \quad c^p=a^4=b^2=a^b*a=(a,c)=c^b*c=1\ \text{[case b in Table 1]}
\end{equation}
for the case where the kernel is $C_4$. 		 

Case a) has an automorphism group isomorphic to
\begin{equation}			
	\text{(entry in table)}  \times  Hol(C_p)
\end{equation}
 or
\begin{equation}
			  C_2 \times  C_2 \times  Hol(C_p)
\end{equation}
and case b) has an automorphism group isomorphic to
\begin{equation}
			       D_4 \times  Hol(C_p),
\end{equation}
where $Hol(C_p)$ is the holomorph of the group $C_p$.\\

	The remaining entries in this table follow the same pattern. The \textquotedblleft entry in table\textquotedblright\ factor will appear in all $C_2$ extensions of $p$-groups by the same order 8 group. 

\section{Groups of order $8p^2$ and their automorphism groups}

	Let us start by considering the case $(C_p \times  C_p) @ D_4$. This is a group of order 72 for $p=3$. In fact, there are several possible groups of this form. Let us for now only concern ourselves with the $C_2$ actions of $D_4$ on the $p$-group. From the above comments we know there are two possible $C_2$ actions that can arise from $D_4$, but each of these actions can act in two different ways on $C_p \times  C_p$, viz:
\begin{gather}
    		    c^p=d^p=(c,d)=a^4=b^2=a^b*a= \notag \\
			(a,c)=(a,d)=(b,c)=d^b*d=1 \qquad \qquad  C_p \times  (C_p @ D_4),
\end{gather}
or
\begin{equation}     
			(a,c)=c^b*c=(a,d)=d^b*d=1  \qquad \qquad (C_p \times  C_p) @ D_4.
\end{equation}
The other case arises from interchanging the actions of $a$ and $b$ on the $p$-group, viz:
\begin{gather}
		   c^p=d^p=(c,d)=a^4=b^2=a^b*a= \notag \\
			(a,c)=(b,c)=d^a*d=(b,d)=1 \qquad \qquad C_p \times  (C_p @ D_4),
\end{gather}
 or
\begin{equation}	
			c^a*c=(b,c)=d^a*d=(b,d)=1 \qquad  \qquad (C_p \times  C_p) @ D_4.
\end{equation}
The automorphism groups for the four different groups of order $8p^2$ are all different but easily specified for all primes $p$:
\begin{gather}
		Hol(C_p) \times  C_{(p-1)} \times  D_4,\\
		Hol(C_p \times  C_p) \times  D_4,\\		
		Hol(C_p) \times  C_{(p-1)} \times  C_2 \times  C_2 ,\\
		Hol(C_p \times  C_p) \times  C_2 \times  C_2.
\end{gather}
In fact, the general case looks like this: 
\begin{equation}
		c_1^p=c_2^p=\cdots =c_n^p=a^4=b^2=a^b*a=\cdots=1
\end{equation}
with $a$ commuting with the generators $c_1,\ldots,c_n$ and $b$ acting as an element of order 2 on the first $m$ generators $c_1,\ldots,c_m$, and commuting with the last $n-m$ generators of order $p$. In this case the automorphism group of this group is
\begin{equation}
		Hol(C_p \times  C_p \times \cdots \times C_p) \times  \textup{GL}((n-m),p) \times  D_4,                
\end{equation}
where we have $m$ $C_p$'s in the holomorph. If the generator $b$ were commuting with the generators $c_1,\ldots,c_n$ and the generator $a$ were acting on the generators $c_1,\ldots,c_n$ instead, then the resulting automorphism group would be obtained by just replacing the $D_4$ factor by the factor $C_2 \times  C_2$ in (15). In a like manner one can find similar patterns in other groups or with other actions, e.g., $C_2 \times  C_2$, $C_4$ or an order 8 group. \\

	The groups of order $8p^2$ that arise by a $C_2$ action of the 2-group on the $p$-group can be easily written out, and their automorphism groups also follow a very simple pattern. Namely, for the automorphism groups, we determine the $C_2$ action on the $p$-group, pull out the corresponding entry in Table 1, and the \textquotedblleft $p$-factor\textquotedblright\ is now either 
\begin{equation}
	      C_{(p-1)} \times  Hol(C_p)    
\end{equation}
if the $C_2$ acts only on one of the $C_p$'s or  
\begin{equation}
		  Hol(C_p \times  C_p)    
\end{equation}
if the $C_2$ acts on both of the $C_p$'s.\\

	When the $p$-group is cyclic, the groups behave just like the case of order $8p$. Namely, just use the \textquotedblleft invariant factor\textquotedblright\ in Table 1 with $Hol(C_{(p^2)})$ and you have the automorphism groups for these cases.\\

	The cases involving an order 4 or higher-order action are listed in Table 2. In the case of $C_2 \times  C_2$ actions the \textquotedblleft $p$-factor\textquotedblright\ is one of two types: $Hol(C_p) \times  Hol(C_p)$, or a group that can be represented as a wreath product \cite{12}. The automorphism groups coming
from a $C_2 \times  C_2$, $C_4$ or a $C_4 \times  C_2$ action are explicitly written out in Table 2. The cases arising from the $D_4$ or $Q_2$ actions depend upon $p$ and are given in the notes to Table 2 (Appendix 1). \\

	Table 2 also gives in an abbreviated form the presentations of the groups of order $8p^2$. To obtain the required presentation from this table a few examples will show how to reconstruct the presentations. Look at the $(C_2 \times  C_2)$ image from $C_4 \times  C_2$, for the 
case [ab,b]. The group whose automorphism group is $[Hol(C_p) \times  C_2] wr C_2$ has the presentation:
\begin{equation}
		a^4=b^2=(a,b)=c^p=d^p=(c,d)=c^a*c=c^b*c=(a,d)=d^b*d=1.
\end{equation}
For the group $(C_4 \times  C_2)$ with an order 8 group action, Table 2 gives [a,ab] (with automorphism group $Hol(C_p) wr C_2)$ giving the presentation:
\begin{equation}
		a^4=b^2=(a,b)=c^p=d^p=(c,d)=c^a*c^x=d^a*d=(b,c)=d^b*d=1,
\end{equation}
where $x^4\equiv 1$ mod($p$) (i.e., $x$ is a fourth root of unity \cite{20}).

	The groups without a normal sylow $p$-subgroup are treated in Table 3.\\

	For the case of $p\equiv 1$ mod(2) but not 1 mod(4) or 1 mod(8) we have only two cases of a $C_4$ and one case of a $C_8$ extension. For $p=7$ we have extensions in which the action is full (i.e., by a group of order 8) involving only the groups $C_8$, $D_4$, and $Q_2$. The automorphism group of $(C_p \times  C_p) @ Q_2$ is always a complete group \cite{13}. The automorphism group of the group with a full $C_8$ action is also a complete group, with order 4704 (for $p=7$). The $D_4$ extension's automorphism group is not complete, but has order 2352 (for $p = 7$). \\

	In connection with the study of \textquotedblleft higher\textquotedblright\ orders, e.g., 2-groups of orders 16, 32, and especially of order 64 the following observation is very useful: if the group associated with the extension is a characteristic subgroup of index 2 for the given 2-group, then the automorphism group arising from this extension is given by \cite{14}
\begin{equation}
		Hol(C_p) \times \text{ Automorphism group(2-group).}\quad      
\end{equation}
In connection with the groups of order 192 with a normal sylow 3-subgroup, this means that about 700 of the automorphism groups of this type in this order are already known, without the need for performing any calculations whatsoever \cite{15}. In fact since these groups have a \textquotedblleft natural extension\textquotedblright\ to groups of order 64$p$, $p$ any odd prime, this whole class of automorphism groups is already known without the need to do any calculations! General results like this considerably reduce computer time in the calculation of automorphism groups, and at the same time provide a link between the structure of a group and its automorphism group. One should also note that this simplification does not extend to the cases of the form $C_q @ X$, where $X$ is a $p$-group ($p$ an odd prime, $q\equiv 1$ mod($p$)) and $X$ acts on $C_q$ as an operator of order $p$. An interesting exercise would be to see what happens both in this case and in the cases for 2-groups where the group associated with the extension is a characteristic subgroup of index greater than 2. 
	
\section{Groups of order $8p^2$. The $C_8$ Extensions}

	The most recent discussion of the groups of order $8p^2$ is that given by Zhang Yuanda \cite{16}. Following Zhang Yuanda \cite{16} the groups of order $8p^2$ with the $p$-group being $C_p \times  C_p$ and with the 2-group being cyclic give rise to the cases listed in Table 4.\\

	All of the numbers in Table 4 refer to the listing of the groups as given in the article by Zhang Yuanda \cite{16}. The relations for the groups in Table 4 are given below.\\ 

      For the purpose of determining the number of groups of order $8p^2$ of the form $(C_p \times  C_p) @ C_8$ for a particular prime $p$, it is more convenient to use the following breakdown:
\linebreak

\begin{tabular}{|l|c|}\hline
\multicolumn{2}{|c|}{ Number of groups of the form $(C_p \times  C_p) @ C_8$} \\ \hline
    & number of groups \\ \cline{2-2}   
prime $p$ &  image of 2-group \\  
  &   $C_2 \quad C_4 \quad C_8$ \\ \hline
	   if $p\equiv 1$ mod(2)  &     2 + 1 + 1     \\ \hline
     if $p\equiv 1$ mod(4)  &      2 + 4 + 1       \\ \hline
	       if $p\equiv 1$ mod(8)  &      2 + 4 + 8    \\ \hline
\end{tabular}
\linebreak

The relations for the $C_4$ action when $p$ is not equal to 1 mod(4) can be written as
\begin{equation}
		a^p=b^p=(a,b)=c^8=a^c * b^{-1}=b^c * a=1.
\end{equation}
The relations that Zhang Yuanda gives for a $C_8$ action when $p \equiv 3$, 5 
and 7 mod(8) are:\\

\begin{tabular}{|l|c|l|c|} \hline
\multicolumn{4}{|c|}{$ a^p=b^p=(a,b)=c^8=a^c*b-1=b^c*b-s*a^{-t}=1 $}\\ \hline
  number in ref. [16]  & relation for $s$ and $t$   &  ($p$,$s$,$t$)  & Z($g$)\\ \hline
                            &                                 &           & \\
  16. for $p\equiv 3$ mod(8)& $ s^2\equiv -2$ mod($p$), $t=1 $&  (11,3,1) &  I  \\ 
                      &                     &               & \\
  7. for $p\equiv 1$ mod(4) &$ s=0, t^2\equiv -1$ mod($p$) &  (13,0,8) &  I  \\
                      &                     &  (17,0,4) &  I  \\
                      &                     &           &     \\ 
 17. for $p\equiv 7$ mod(8) &$ s^2\equiv 2$ mod($p$), $t=-1$&  (7,3,-1) &  I  \\
                      &                     &  (31,8,-1) &  I  \\ \hline
\end{tabular}
\linebreak

	Explicit relations for the first few cases of $(C_p \times  C_p) @ C_8$ are given by: 
\begin{equation}
		a^p=b^p=(a,b)=c^8=a^c*a^w*b^x=b^c*a^y*b^z=1,
\end{equation}
where the entries in Table 5 are ($w$,$x$;$y$,$z$). The automorphism groups for these groups are given in Table 6. \\

	The following two sets of presentations yield groups that are not distinguishable by their class structure and automorphism groups alone: \\

\begin{tabular}{|c|c|c|} \hline
\multicolumn{3}{|c|}{    [-2,2], [-2,8], [-2,9],$\longrightarrow  $\#'s 7, 12, 13 in Table 6}\\ \hline
  number of elements & order of elements & number of classes \\ \hline
       289           &        2          &         1         \\ \hline
	  578          &        4          &         2         \\ \hline
	 1156           &        8          &         4         \\ \hline
	  288           &       17          &        36        \\ \hline
\end{tabular}
\linebreak

Here we are using the notation [-2,a] = (-2,0:0,a). \\

\begin{tabular}{|c|c|c|} \hline
\multicolumn{3}{|c|}{[-2,4], [-2,13]$\longrightarrow $ \#'s 11 and 10 in Table 6}\\ \hline
  number of elements & order of elements & number of classes \\ \hline
          17         &        2          &         1        \\ \hline
	    578         &        4          &         2         \\ \hline
	   1156         &        8          &         4         \\ \hline
	    288         &       17          &         38       \\ \hline
	    272         &      34          &        4         \\ \hline
\end{tabular}
\linebreak

This poses the interesting question of how to determine when two groups with different presentations might be isomorphic. Clearly in many 
cases one can compare the conjugacy classes of the two groups to show that they are not isomorphic. In other cases when this does not distinguish between two groups the automorphism groups of these groups may be different. Within each of the sets above neither of these methods enables one to be able to distinguish between these nonisomorphic groups. (An interesting exercise for the reader would be to develop a computational method which would enable one to 
distinguish between these different groups \cite{17}.) Such methods obviously exist. It would be helpful to be able to list the group invariant that enables one to distinguish between these groups. Group theory programming systems such as GAP have a routine that enables its users to determine whether two groups are isomorphic or not. \\

\section{Conclusions}

	We were told by some group theorists when we started this work around 1980 that the study of the automorphism groups would be very difficult since there were very few general results known \cite{18}. We felt at that time if one could determine general presentations for the groups of a particular type (e.g.,of order $8p$), then one should probably be able to determine their automorphism groups with equal precision. As this work shows, we have in fact been able to 
obtain explicit expressions for the automorphism groups for the groups of orders $8p$  and $8p^2$, for all odd primes $p$ (modulo a representation theory problem for the $p\equiv 7$ mod(8) cases). In this work the \textquotedblleft factor\textquotedblright\ that is independent of \textquotedblleft $p$\textquotedblright\ in Aut($G$) appears to be a new group invariant associated with the group $G$ (or Aut($G$)). More details on this aspect of the work will be presented in a paper on the automorphism groups for the groups of order 32$p$.\\

	The situation arising with the matrix representations for the group $<2,3,4>$ discussed in Appendix 2 also shows the limitations of a purely computational approach to the subject. If one 
had restricted the calculations to primes $p < 103$, then one would have been led to believe that one could always find a matrix representation for the group $<2,3,4>$ using matrices of the form (\ref{A21}) in Appendix 2. As is apparent from the calculations mentioned in Appendix 2 this is not the case. \\

	The use of computers in mathematical research is widespread. The real value of a computer in many areas of mathematics, like group theory among others, seems to us to be in providing a large number of examples, or conjectures if you will, to enable the mathematician to formulate new theorems which may be proven by other noncomputational means. In the case of group theory this is certainly a desirable goal. The computer like all finite machines has its limits. In some of the computer runs to be reported in other papers the time taken to compute one automorphism group may be as much as 24 hours. Clearly as one goes to larger and larger orders other less brute force methods become more and more relevant. The insight provided by the lower-order calculations should lead one to a better understanding of just what is likely to be true and what is false. A single example can show that a \textquotedblleft theorem\textquotedblright\ is false, but no finite
number of cases can show that it is true. Hence the intelligent interactions of the mathematician and the computer hopefully provide a guide to the formulation and solution of otherwise difficult problems. The authors' work on the computation of the automorphism groups is an attempt to provide the first step in showing just how useful the computation of automorphism groups can be in providing a guide as to what the general structure of these automorphism groups is. The author hopes the conjectures or \textquotedblleft theorems\textquotedblright\ proposed above show just how useful this interaction of computers with the more standard methods of group theory can be for future research.\\

\pagebreak

\section{Appendix 1 \cite{20}. Notes for Table 2}

\underline{	Mainly presentations for automorphism groups of $D_4$ and $Q_2$ images.}\\

	 	The group [144] in Table 2, for $p > 3$ means 
\begin{equation*}			
Aut[ (C_p \times  C_p) @ C_4] 
\end{equation*}
	when  $p\equiv 3$, 7 mod(8), or $H(p^2)$. Here $H(p^n)$ is the group of all mappings $x \mapsto  a*x^t + b$ where $a$ (different from zero) and $b$ are elements of the Galois field, GF$(p^n)$, and $t$ is an element of the Galois group of GF$(p^n)$. Here $p$ is a prime number. The order of the group $H(p^n)$ is $p^n * (p^n - 1) * n$.  (See D. Robinson \cite{19}.)\\

\underline{	Notes for $C_4$ column:}\\

	For $p\equiv 1$ mod(4) we have the following four cases for a $C_4$
	action on the $p$-group, and the [144] factor being replaced by 
	the factors below arising from the four nonisomorphic cases, 
	which in lowest order ($p = 5$, and with the $C_4$ action in $C_4 \times  C_2$) 
	take the form:\\

\begin{tabular}{|c|c|c|} \hline
\multicolumn{2}{|c|}{ Group ($p =  5$ case) } & Automorphism groups\\ \cline{1-2}
      Group     & relations ($p = 5$)   &          for $p > 5$       \\ 
                  &  $ a^5=b^5=(a,b)=c^4= $   &                        \\ \hline
  $ (C_5 @ C_4) \times  C_5 $&$   a^c*a^2=(b,c)=1;  $&$  Aut(g) = Hol(C_p) \times  C_{(p-1)}$ \\
      $Hol(C_5) \times  C_5 $&           &                           \\ 
                     &                       &                           \\
     $(C_5 \times  C_5) @ C_4  $&$  a^c*a^2=b^c*b=1;     $&$ Aut(g) = Hol(C_p) \times  Hol(C_p)$\\
                     &                       &                           \\
     $(C_5 \times  C_5) @ C_4  $&$  a^c*a^2=b^c*b^2=1;   $&$ Aut(g) = Hol(C_p \times  C_p)$     \\
                     &                       &                           \\
    $ (C_5 \times  C_5) @ C_4  $&$  a^c*a^2=b^c*b^3=1;   $&$ Aut(g) = Hol(C_p) wr C_2$.   \\
                    &or $ a^c*(b^{-1})=b^c*a=1$; &                          \\ \hline
\end{tabular}
\linebreak

	       The last two groups have the same class (order) structure  
	and so cannot be distinguished from their conjugacy classes alone.
	As can be seen above though, they can be distinguished from the
	fact that they have different automorphism groups.\\

		The sequence of groups of order $4p^2$ with the presentations
\begin{equation}
		    a^p=b^p=(a,b)=c^4=a^c*(b^{-1})=b^c*a=1
\end{equation}
	have complete groups for their automorphism groups for 
	$p$ =(3,5,7,11). For primes $p < 23$ the automorphism groups 
	and their orders are:\\

\begin{tabular}{|c|c|c|} \hline
	        $p$      &     order of group        &automorphism group   \\ \hline
	         3      &$     144 = 2^4 * 3^2      $&$     H(3^2)      $    \\ \hline
	         5      &$     800 = 2^5 * 5^2      $ &$  Hol(C_5) wr C_2$      \\ \hline
	         7      &$    4704 = 2^5*3*7^2      $ &$     H(7^2)      $    \\ \hline
	        11      &$  29,040 = 2^4*3*5*11^2   $ &$     H(11^2)     $    \\ \hline
	        13      &$  48,672 = 2^5*3^2*13^2   $ &$  Hol(C_{13}) wr C_2 $  \\ \hline
	        17      &$ 147,968 = 2^9 * 17^2     $ &$  Hol(C_{17}) wr C_2 $    \\ \hline
	        19      &$ 259,920 = 2^4*3^2*5*19^2 $ &$   H(19^2)        $   \\ \hline
\end{tabular}
\linebreak

		 If $p\equiv 1$ mod(4), the automorphism group of the group whose 
	presentation is
\begin{equation}
			a^p=b^p=(a,b)=c^4=a^c*(b^{-1})=b^c*a=1
\end{equation}
is $Hol(C_p) wr C_2$; otherwise, it is $H(p^2)$ of order $p^2 * (p^2 - 1) * 2$. In either case,
Aut($g$) is complete.
\begin{quotation}		
	Note. Not all products of the form $(Hol(C_p) wr C_2)$ are complete
	groups. The group $Hol(C_3) wr C_2$ has order 72 and is not complete, 
 	but $Hol(C_7) wr C_2$ is a complete group.
\end{quotation}
	Notes for order 8 column:\\

		\underline{a). For the $C_8$ cases see Table 4.}\\

		\underline{b). $D_4$ cases.}\\
			 The presentation for this set of groups is:
\begin{gather*}
			a^4=b^2=a^b*a=c^p=d^p=(c,d)=\\
				c^a*d=d^a*(c^{-1})=c^b*c=(b,d)=1.
\end{gather*}
		   The automorphism groups for this sequence of groups
		have the following orders:
\linebreak

\begin{tabular}{|c|c|c|} \hline
		   $p$     &  order Aut($g$)   &      comments         \\ \hline
		    3     &      144        &  complete group       \\ \hline
	          5     &      800        &  complete group       \\ \hline
		    7     &    2,352        &  not complete   (1)  \\ \hline
		    11    &    9,680        &  complete group       \\ \hline
		    13    &   16,228        &  complete group       \\ \hline
		    17    &   36,992        &  not complete    (2)  \\ \hline
		    19    &   51,984        &  complete             \\ \hline
		    23    &   93,104        &  not complete    (3)  \\ \hline 
		  \multicolumn{3}{|l|}{(1) Aut(2,352) has order 4,704. }\\
		   \multicolumn{3}{|l|}{}\\
		 \multicolumn{3}{|l|}{ (2) $36,992\rightarrow 73,984\rightarrow Hol(C_{17}) wr C_2$ (complete).}\\
		   \multicolumn{3}{|l|}{}\\
		 \multicolumn{3}{|l|}{(3) Aut here has order 186,208.}  \\ \hline
\end{tabular}
\linebreak

		The orders of these Aut($g$)'s follow the pattern $8*p^2*(p-1)$, the $p\equiv 3$ or 5 mod (8) being complete. If  $p\equiv 1$ or 7 mod(8), then the automorphism group is not a complete group.
The relations for the $p$ = 3, 5, 7, and 1  mod 8 are:
\linebreak

\begin{tabular}{|c|l|c|} \hline
\multicolumn{3}{|c|}{     Table  D   } \\ \hline
\multicolumn{3}{|c|}{  automorphism groups for $(C_p \times  C_p) @ D_4$ }\\ \hline             
        primes   &        presentations             &modular relations \\ \hline
  && \\
          3 mod(8)  &$ a^8=b^2=a^b*a^5=c^p=d^p=(c,d)=$   &$16*x^8 \equiv 1$ mod($p$)\\
                    &$ c^a*(c^-x)*d^x=d^a*(c^{-x})*(d^{-x})= $&  $y=(p-1)/2$    \\
        $ (C_p \times  C_p)  $&$ (b,c)=d^b*d  =  f^y=(a,f)= $&$  t^y \equiv 1$ mod($p$)\\
        $@ (QD_8 \times  C_y)$&$ (b,f)=c^f*(c^{-t})=d^f*(d^{-t})=1 $&\\ \cline{3-3}
                    &                                   &  ($p,x,y,t$)      \\
                    &                                   &    (3,1,1,1)    \\
                    &                                   &   (11,4,5,3)    \\
                    &                                   &   (19,3,9,4)    \\ \hline
     &&\\
           5 mod(8)   &$   a^p=a^b*(a^{-x})=b^4=c^2=   $& $x^4 \equiv 1$ mod($p$)  \\
                      &$   (a*c)^2*((a^{-1})*c)^2=        $    &  $q=(p-1)/4$       \\
         $[(C_p @ C_4)$   &$  a*c*b*c*(a^{-1})*c*(b^{-1})*c=     $   & $ y^q \equiv 1$ mod($p$)  \\
          $wr C_2] @ C_q$ &$   (b*c)^2*((b^{-1})*c)^2=        $  & \\ \cline{3-3}
                      &$  d^q=a^d*(a^{-t})=b^d*a*(b^{-1})*a=$    &    ($p,x,q,y$)     \\
                      &$  (c,d)=1                        $  &   (13,-5,3,3)    \\  
                      &                                    &   (29,12,7,-13)  \\ \hline
                 && \\ 
          7 mod(8)   &$ a^8=b^2=a^b*a=c^p=d^p=(c,d)=      $ &$16*x^8 \equiv 1$ mod($p$) \\
                      &$ c^a*(c^{-x})*d^x=d^a*(c^{-x})*(d^{-x})=$  &  $y=(p-1)/2$      \\
          $(C_p \times  C_p)$   & $(b,c)=d^b*d  =  f^y=(a,f)=      $   & $ t^y \equiv 1$ mod($p$) \\
         $@ (D_8 \times  C_y)$  & $(b,f)=c^f*(c^{-t})=d^f*(d^{-t})=1  $ & \\ \cline{3-3}
                      &                                    &  ($p,x,y,t$)       \\
                      &                                    &    (7,2,3,4)     \\
                      &                                    &   (23,9,11,4)    \\ \hline       
                  && \\
          1 mod(8) *  & $a^p=b^q=a^b*(a^{-x})=c^4=d^2=(a,c) $ &  $q=(p-1)$         \\
	 	      & $=(b,c)=(b,d)=(a*d)^2*((a^{-1})*d)^2 $ & $ x^q \equiv 1$ mod($p$)  \\
	              &$= b^t*(d^{-1})*(c^{-1})*(d^{-1})*(c^{-1})=1 $ &   $t = q/4$       \\
                      &                                    &                  \\ \hline
\multicolumn{3}{|c|}{ * conjecture, see text below.}        \\ \hline
\end{tabular}
\linebreak \\

 		If we replace the term $a^b*a^5$ in the $p \equiv 3$ mod(8) presentation and the
	term $a^b*a$ in the $p \equiv 7$ mod(8) presentation above with $a^b*a^s$ with
	$s$ given by 
\begin{equation*}
			p \equiv (8-s)\,\, \textup{mod}(8),
\end{equation*}
	we then have just one form that is valid for both the $p \equiv 3$ mod(8) 
	and the $p \equiv 7$ mod(8) presentations. \\

		A four-generator permutation representation of this group 
	for $p$ = 13, of degree 26 is:
\begin{align*}
			a&=(1,2,3,4,5,6,7,8,9,10,11,12,13);\\
			b&=(1,8,12,5)(2,3,11,10)(4,6,9,7);\\
			c&=(1,14)(2,15)(3,16)(4,17)(5,18)(6,19)(7,20)\\
			  &\qquad	(8,21)(9,22)(10,23)(11,24)(12,25)(13,26);\\
			d&=(2,10,4)(3,6,7)(5,11,13)(8,12,9)\\
  			  &\qquad  (15,23,17)(16,19,20)(18,24,26)(21,25,22).\\
\end{align*}
		For the $p \equiv 1$ mod(8) case we have only been able to verify 
	the presentation for the first case, namely that for $p$ = 17,
\begin{gather*}
			a^{17}=b^{16}=a^b*a^3=c^4=d^2=(a,c)=(b,c)=(b,d)=\\
			(a*d)^2*((a^-1)*d)^2=b^4*((d^-1)*(c^-1))^2=1.
\end{gather*}
	We conjecture that the general matrix representation of the group
	of order $8(p-1)$ is given by
\begin{equation}
b= \left(\begin{matrix}
       z & 0 \\
       0 & z \end{matrix} \right),
    \quad 
c = \left( \begin{matrix}
        1 & 0 \\
        0 & x \end{matrix} \right), \quad
d = \left(\begin{matrix}
      0 & 1 \\
      1 & 0 \end{matrix} \right),
\end{equation}
	where
\begin{gather*}    
				x^4 \equiv 1\,\, \textup{mod}(p), \\
				z^{(p-1)} \equiv 1\,\, \textup{mod}(p).
\end{gather*}
	The matrix $b$ is the matrix representation for the center of the 
	group GL(2,$p$) (see [20]), and the group generated by $<c,d>$ is 
	$C_4 wr C_2$. The group $<b,c,d>$ appears to be the same group arising in 
	the automorphism groups for the groups of order $16p^2$ with a $D_4$ 
	or $Q_2$ action on the group $C_{17} \times  C_{17}$. \\
		
\underline{	c). $Q_2$ image case.}\\

		The presentations of these groups of order $8p^2$ can
	be read from the representation of $Q_2$ in GL(2,$p$).\\

		A matrix representation of $Q_2$ for the presentation
\begin{equation}		
			a^4=b^4=a^2*b^2=a^b*a=1
\end{equation}
is
\begin{equation}
a = \left(\begin{matrix}
      x & y \\
     y & -x  
      \end{matrix} \right),  \qquad 
   b = \left( \begin{matrix}
         0 & 1 \\
         -1 & 0 
     \end{matrix} \right),
\end{equation}	    
	where $(x^2 + y^2) \equiv -1$ mod($p$).\\

		The automorphism groups of these groups have the structure [21]:\\
\begin{equation*}
			(C_p \times  C_p) @ \text{(group of order}\,\, 24(p-1)\,\,).\\
\end{equation*}
	The group of order $24(p-1)$, in this sequence of automorphism groups, 
	depends upon the prime $p$ as follows:\\

\begin{tabular}{|l|cl|l|} \hline
  prime     & \multicolumn{2}{|c|}{ quotient group}&   \\ \hline
	       $p \equiv 3$ mod(8)&   GL(2,3)   &$ \times  C_q$ & $q=(p-1)/2$   \\
	       $p \equiv 5$ mod(8)&   SL(2,3)$ @ C_4 $&$  \times  C_q $& $q=(p-1)/4$   \\
	       $p \equiv 7$ mod(8)&   $<2,3,4>$& $\times  C_q $& $q=(p-1)/2$  \\
	       $p \equiv 1$ mod(8)& GL(2,3)$ @ C_{(2^n)}$&$ \times  C_q $& $(p-1)=q*2^n$\\ \hline
\end{tabular}
\linebreak

		In the cases where  $p \equiv 1,$ 3, or 5  mod(8) these groups of order $24(p-1)$ can be written in the form:
\begin{equation*}
			\textup{SL}(2,3) @ C(p-1).
\end{equation*}
with the presentation:  
\begin{gather}
			a^2=b^3=(a*b*a*b*a*(b^{-1}))^2=(a,c)=(b,c)=  \notag \\
				c^x=c^y*(a*b)^4=1,                        
\end{gather}
where $x=(p-1)$ and $y=x/2$. A matrix representation for this series of groups is given by:
\begin{equation}		    
a = \left(\begin{matrix}
         0 & s \\
         t & 0 \end{matrix} \right), \qquad
b = \left( \begin{matrix}
     -1 & 1 \\
     -1 & 0 \end{matrix} \right), \qquad
c = \left( \begin{matrix}
	z & 0\\
      0 & z \end{matrix} \right), 
\end{equation}
where the matrix $c$ is just the center of the group GL(2,$p$) and $H =$ $<a,b>$ is a representation of GL(2,3) by $2 \times  2$ matrices over the field GF($p$). For the $p \equiv 1$ mod(8) case we have for the first few cases:
\begin{align}
			&p = (17,41,73,89,97,113);\\
			&s = (2,3,10, \ldots);\\
			&t = (9,14,22, \ldots);\\
			&z = (3,6,5,3,5,3, \ldots). 
\end{align}
In the general case we have $s$, $t$ and $z$ being given by the following relations:
\begin{gather}
			s^8 \equiv 1\,\, \textup{mod}(p),\\
			t = (p+1)/s,
\end{gather}
and
\begin{equation}
			z^{(p-1)} \equiv 1\,\, \textup{mod}(p)
\end{equation}
(see [20]). In the other case, the group GL(2,3) is replaced by the group $<2,3,4>$. In this case we have the following presentation:
\begin{gather}
	(a^{-2})*b^3=(a^{-2})*c^4=(a^{-1})*b*c= \notag \\
		        d^x=d^y*a^2=(a,d)=(b,d)=(c,d)=1.
\end{gather}
	Explicit presentations for the $Q_2$ image cases:\\
\begin{quotation}
			Nota bene in the presentations below that the roles of
				$a$ and $b$ are interchanged in the 
				$p \equiv 3$ and 5 mod(8) cases.\\
\end{quotation}

\begin{tabular}{|c|l|c|} \hline
\multicolumn{3}{|c|}{     Table  Q   } \\ \hline
\multicolumn{3}{|c|}{  automorphism groups for $(C_p \times  C_p) @ Q_2$ }\\ \hline             
        primes   &        presentations             &modular relations \\ \hline
  && \\
         3 mod(8)   &$ a^3=b^2=((a*b)^2*(a^{-1})*b)^2= $ & $q = (p-1)/2$      \\
    	              &$ c^p=d^p=(c,d)=c^a*c*(d^{-1})=d^a*c=$ &$x^{(p-1)}\equiv 1$ mod($p$) \\
                    &$ c^b*(c^{-1})*(d^{-x})=d^b*d=          $& $ y^q \equiv 1$ mod($p$)  \\
                    &$ e^q=c^e*(c^{-y})=d^e*(d^{-y})=       $ &                  \\
                    &$  (a,e)=(b,e)=1                  $ &                  \\ \hline
    && \\
	   5 mod(8)  &$ a^3=(a,b^2)=(a*b)^4=(a*(b^{-1}))^4=$ &   $q=(p-1)/4$      \\
                    &$ c^p=d^p=(c,d)=c^a*c*(d^{-1})=d^a*c= $& $ x^4 \equiv 1$ mod($p$)  \\
                    &$ c^b*(d^{-x})=d^b*(c^{-1})=            $& $ y^q \equiv 1$ mod($p$) \\
                    &$     e^q=(a,e)=(b,e)=c^e*(c^{-y})=   $&                  \\
                    &$         d^e*(d^{-y})=1             $&                  \\ \hline
	 7 mod(8)   &    see discussion below           &                  \\ \hline
       && \\
         1 mod(8)   &$ a^2=b^3=(a*b*a*b*a*(b^{-1}))^2=  $&   $q=(p-1)$        \\ 
                    &$  (a,c)=(b,c)=c^q=c^t*(a*b)^4=    $&   $t=q/2$          \\
                    &$   d^p=e^p=(d,e)=d^a*(e^{-x})=    $ &                  \\
                    &$   e^a*(d^{-y})=d^b*d*(e^{-1})=e^b*d= $&$ x^8 \equiv 1$ mod($p$) \\
                    &$      d^c*(d^{-z})=e^c*(e^{-z})=1  $& $y = (p+1)/s$      \\
                    &                                   & $z^{(p-1)}\equiv 1$ mod($p$) \\
                    &                                   &                  \\ \hline
\end{tabular}
\linebreak

	The case of $p \equiv 7$ mod(8). The automorphism groups in this sequence are given by the following  presentation:
\begin{gather} 
			a^3=b^4=(a,b^2)=a*b*(a*(b^{-1}))^3= \notag  \\
			c^p=d^p=(c,d)=c^a*c*(d^{-1})=d^a*c= \notag \\
			c^b*c^v*d^x=d^b*c^y=d^w=  \label{App} \\
			e^q=(a,e)=(b,e)=c^e*(c^{-z})=d^e*(d^{-z})=1, \notag
\end{gather}
where
\begin{gather*}
(p,v,w,x,y)=(7,1.-1,1,5),\quad (23,1,-1,7,3),\\
		(31,1,-1,12,5),\quad (47,1,-1,18,26),\quad (71,1,-1,7,20).....
\end{gather*}
  and $q=(p-1)/2$, $z^q \equiv 1$ mod($p$).\\

		Here the matrices are
\begin{equation}\label{A12}
a = \left(\begin{matrix}
      -1 & 1 \\
      -1 & 0 \end{matrix} \right), \qquad
b = \left(\begin{matrix}
      v & x \\
      y & w \end{matrix} \right).
\end{equation}
Here $a$ and $b$ generate the group $<2,3,4>$. \\

		An alternate representation for this automorphism group is 
	based upon the group $<2,3,4> \times  C_q$ and is generated by matrices of the 
	form (\ref{A12}), and is (for $p$ = 7)
 \begin{gather}\label{A13a}
		   a^3=b^{(4q)}=(a,b^2)=(a*(b^{-1}))^4*(b^{-(p-5)})=\notag \\
		   c^p=d^p=(c,d)=c^a*c*(d^{-1})=d^a*c=  \\
		   c^b*((c^{-1})*(d^{-x})=d^b*d*(c^{-y})= \notag \\                                             
	           (a*b)^4*b^2=1,  \notag 
\end{gather}
and for $p > 7$ we have
\begin{gather}\label{A13b}
         	 a^3=b^{(4q)}=(a,b^2)=(a*(b^{-1}))^4*(b^{-(p-5)})= \notag \\
		   c^p=d^p=(c,d)=c^a*c*(d^{-1})=d^a*c=    \\
		   c^b*((c^{-1})*(d^{-x})=d^b*d*(c^{-y})= \notag \\                                                                       
		   (a*b)^2*(a^{-1})*(b^{-1})*(a*(b^{-1}))^2*(a^{-1})*b=1. \notag
\end{gather}
	In this form [b] has order $4*q=4*(p-1)/2$. For many (but not 
	all!) primes of the form $p \equiv 7$ mod(8) one can represent the 
	generator [b] for the presentation (\ref{App}) with $v=-w=1$, and  
	$x*y \equiv -2$ mod($p$). One can also find representations for [b] for the 
	form (\ref{A23}) with $v=-w=1$. I do not have closed form expressions for the 
	values of $x$ and $y$, for either form of the above presentations.\\
  
		This is a problem in the representation theory of $<2,3,4>$ in 
	terms of $2 \times  2$ matrices over the field GF($p$). See Appendix 2 for 
	details. \\

\section{Appendix 2} 

\underline{Question on Matrix Representations of
		the Coxeter Group $<2,3,4>$ by $2 \times 2$}\\
 \underline{Matrices over GF($p$).}\\

	The Coxeter group $<2,3,4>$ arises in connection with  
	getting the presentations for the automorphism groups of
	the groups
\begin{equation}
		(C_p \times  C_p) @ Q_2 \qquad \text{when $p \equiv 7$ mod(8)}.
\end{equation}
	For many cases the following matrices will give a matrix 
	representation for the Coxeter group $<2,3,4>$:
\begin{equation}\label{A21}
a=\left(\begin{matrix}
              -1 & 1 \\
              -1 & 0 
\end{matrix} \right), \qquad
 b=\left(\begin{matrix}
              1 & x \\
              y & -1 
\end{matrix} \right). 
\end{equation}
Here the matrix $a$ has order 3 and $b$ is of order 4. The presentation for this matrix representation is 
\begin{equation}\label{A22}
		a^3=b^4=(a,b^2)=a*b*(a*(b^{-1}))^3=1      
\end{equation}
when the entries in the matrix $b$ obey the relation:
\begin{equation}
			x*y \equiv -2\,\, \textup{mod}(p).
\end{equation}
	Problem. Find an algebraic method for determining the values
	of $x$ and $y$ in the above matrix $b$. The following are the values
	found by trial and error (actually a computer run for various
	primes).\\

\begin{tabular}{|c|c|c|c|c|c|} \hline
prime & ($x$,$y$) & prime & ($x$,$y$) & prime & ($x$,$y$) \\ \hline
            7   &   (1,5)  &   103   &   none   &   223  &   (57,43)  \\
                &   (2,6)  &         &          &        & (101,117)  \\
           23   &   (7,3)  &   127   &   none   &        & (106,122)  \\
                &  (20,16) &         &          &        & (180,166)  \\
           31   &   (12,5) &   151   &   none   &   239  &  (31,131)  \\
                &  (26,11) &         &          &        &  (46,187)  \\
           47   &  (18,26) &   167   & (54,68)  &        &  (52,193)  \\
                &  (21,29) &         & (99,113) &        & (108,208)  \\
                &  (20,14) &   191   &  (8,143) &   263  &    none    \\
                &  (33,27) &         & (42,100) &        &            \\
           71   &   (7,20) &         & (48,183) &   271  &  (114,19)  \\
                &  (51,64) &         & (91,149) &        &  (209,35)  \\
           79   &  (19,29) &   199   &  (12,33) &        &  (236,62)  \\
                &  (50,60) &         & (166,187)&        &  (252,157) \\ \hline
\end{tabular}
\linebreak

	Matrix representations of the form (\ref{A21}) do not exist for the 
	primes 103, 127, 151 and 263. The values listed in the above
	table are apparently the only possible values for $x$ and $y$ 
	yielding the above representation. \\

		Alternate forms are required for the primes 103, 127,
	151, and 263 and probably others larger than 263. A few 
	possible choices for the matrix $b$, which together with $a$
	above obeys the relations (\ref{A22}), are:
\begin{equation}
b = \left( \begin{matrix}
       v & x \\
       y & w 
    \end{matrix} \right),
\end{equation}
	where $(v,x,y,w)$ are given in the following table:\\

\begin{tabular}{|c|c|} \hline
          prime   &              ($v$,$x$,$y$,$w$) values                  \\ \hline
            103   &   (99,99,30,4), (43,43,48,60), (62,21,18,41)...\\
            127   &  (65,65,19,62), (78,29,123,49) (53,117,27,74)..\\
            151   & (133,79,15,18), (77,12,135,74), (52,15,21,99)..\\
            263   & (82,82,11,181), (164,28,82,99), (191,16,54,72).\\ \hline
\end{tabular}
\linebreak

	The values $(v,x,y,w)$ given are just a sample; there are a large 
	number of quadruples that will work here. They were found by 
	using the matrix representation for $<2,3,4> \times  C_q$ given below and
	then raising the matrix $b$ to the power $q$.\\

		For the primes 103, 127, 151, and 263 (as well as for the
	other primes $p$ = 23, 31, \ldots), alternate matrix representations 
	were found for the group $<2,3,4> \times  C_q$. A selected sample of the 
	values for $(x,y)$ that appear in the matrix $b$ (which now has order 
	$4*(p-1)/2$) follows:
\begin{equation}
b = \left( \begin{matrix}
      1 & x \\
      y & -1
      \end{matrix}
      \right),   
        \qquad \quad \text{where $b$ has order 4$\left(\frac{p-1}{2}\right)$.}
\end{equation}
	Some choices for $(x,y)$ are given in the following table:\\

\begin{tabular}{|c|c|c|c|}\hline
	        prime &   ($x$,$y$)  &  prime  &   ($x$,$y$)   \\ \hline
	        7 [4] &   (1,4)  &  79 [1] &  (15,22)   \\
        	        &   (3,3)  &   103 * &  (1,44)    \\
        	        &   (3,6)  &         & (100,44)   \\
        	        &   (4,4)  &         &            \\
        	 23 [22]&   (1,9)  &   127 \* &  (1,56)    \\
        	        &   (1,16) &         &  (2,26)    \\
        	        &   (2,7)  &         &  (24,2)    \\
        	 31 [16]&   (3,14) &   151 \* &  (4,125)   \\
        	        &   (4,26) &         &  (25,95)   \\
        	 47 [23]&   (1,19) &         &  (52,193)  \\
        	    \*   &   (1.39) &         & (108,208)  \\
        	        &  (12,22) &  263 *  &   (1,29)   \\
       	        &   (17,3) &         & (29 ,200)  \\
        	 71 [8] &  (17,13) &         &            \\
        	        &  (61,43) &         &            \\ \hline
 \multicolumn{4}{|c|}{The numbers in the square brackets indicate the} \\
 \multicolumn{4}{|c|}{number of solutions found. A * means the} \\
 \multicolumn{4}{|c|}{run was truncated before all solutions}   \\
 \multicolumn{4}{|c|}{were found. In the case of $p$ = 103 all were}  \\
 \multicolumn{4}{|c|}{found but it was a very large number} \\
 \multicolumn{4}{|c|}{and was not counted. No computer runs }  \\
 \multicolumn{4}{|c|}{were done for the other primes.    }      \\ \hline
 \end{tabular}
\linebreak\\

		The presentations for the groups $<2,3,4> \times  C_q$ are given by (for $p$ = 7):
\begin{equation}
		a^3=(a,b^2)=(a*b)^4*b^2=(a*(b^{-1}))^4*b^{-2}\quad  \text{for  $p$ = 7,}
\end{equation}
and  for $p > 7$,
\begin{equation}\label{A23}
		a^3=(a,b^2)=(a*(b^{-1}))^4*(b^{-x})=  \\                     
		(a*b)^2*(a^{-1})*(b^{-1})*(a*(b^{-1}))^2*(a^{-1})*b=1.
\end{equation}
Here $x = p-5$. These presentations for the group $<2,3,4> \times  C_q$ have been checked out for primes up to 151 and for $p$ = 263. \\

	This thus poses the interesting question in the theory of group representations for the group $<2,3,4>$: 
\begin{quotation}
		 What algebraic relationship does the pair $(x,y)$ 
		obey such that the presentation (\ref{A22}) (or \ref{A23}) is 
		satisfied? More generally, what conditions on the
		elements of the following order 4 matrix
\begin{equation}
b' = \left( \begin{matrix}
       v& x \\
       y & w \end{matrix} \right)
\end{equation}			
are required in order for $<a,b'>$ to obey the presentation or relations (\ref{A22})?\\

If one has access to a programming system such as Maple or Mathematica one may be able investigate this problem rather easily; otherwise, it could be a fairly messy algebra problem. Remember what you are looking for is a representation valid for all primes $p \equiv 7$ mod(8).
\end{quotation}

\section{Acknowledgements}

	This work was started at Michigan State University around 1980. Additional computations were carried out during the authors' stay at the University of Rhode Island and during several summers at Syracuse University's High Energy Particle Groups' Vax cluster, and was completed at Brown University with the Department of Linguistics and Cognitive Sciences DEC computers. The order 72 groups were done at Michigan State University, the orders 200 and 392 were done at 
URI and in Syracuse. The higher-order cases were done at Brown University.  \\

	We would like to thank Drs. M. Goldberg and G. C. Moneti for making available time on the Syracuse University DEC cluster, and Mr. Carl Brown and later Ms. Judith Reed for assistance in getting CAYLEY running there. At Brown University, Dr. James Anderson has allowed us to use a DEC 6000-510 to finish the work reported here as well as a great deal of additional work that will be reported on in other papers. Ms. Margaret Doll at Brown University has been extremely helpful in getting the various versions of CAYLEY up and running on the DEC computers at Brown University and assisting in clearing up other computer problems as they arose in the course of this work. \\

	A special expression of thanks must go to Dr. John Cannon for making available to us the programming system CAYLEY, without which none of this work could have been done. We very much regret that we (Dr. Cannon and the authors) live so very far away from each other. Dr. Cannon is a very fine gentleman and it would have been a great pleasure to have been able to share the results (and problems found in this work using CAYLEY) with him on a more personal basis.

\pagebreak

\begin{tabular}{|c|cc|cc|c|c|}\hline
\multicolumn{7}{|c|}{    Table 1            }\\
\multicolumn{7}{|c|}{Groups of order $8p$  and their automorphism group factors }\\ \hline
  	     &  \multicolumn{5}{|c|}{Normal sylow $p$ cases; }&  Normal sylow  \\
          & \multicolumn{5}{|c|} {image of 2-group  }& 2-subgroup cases\\ \cline{2-6}
2-group & \multicolumn{2}{|c|}{$C_2$}&\multicolumn{2}{|c|}{$C_4$}&$C_8$& \\ \hline
 $ C_8$  & [a]&$ C_2 \times  C_2$ & [a]& $  C_2$  &    I   &               \\ \hline
 $ C_4 \times  C_2$& [a] & $D_4$      & [a] &  $C_2$  &         &            \\
  	     & &             &    &       &         &                     \\
  	     & [b]&  $C_2$       &   &        &         &              \\ \hline
  $C_2 \times  C_2$ & [a]&  $S_4$      &  &       &       & yes,  $G=A_4 \times  C_2$       \\
    $\times  C_2$ &  &            &   &        &    	&    Aut($G$) = $S_4$        \\
  	     & &             & &          &         &                       \\
  	     & &             & &          &         & yes, Frobenius group  \\
 	     & &             & &          &         &     of order 56       \\
  	     & &             & &          &         &  $|$Aut($G$)$|$ = 168*      \\ \hline
     $D_4$   & [a]& $C_2 \times  C_2$ &  &         &         &                       \\
  	     &   &           &&         &       &                     \\
  	     & [b]&  $D_4$      &&           &         &                       \\ \hline
     $Q_2$   & [b] & $D_4$      & &          &       & yes, $G=$ SL(2,3)        \\
  	     &    &          & &          &         &    Aut($G$) = $S_4$        \\ \hline
 \end{tabular}

\begin{tabular}{|c|} \hline
                       Notes for Table 1                             \\ \hline
 	* This order 168 group is the complete group of this order.   \\ 
    This group is isomorphic to a degree 8 permutation group         \\
   with generators:                                                   \\
           $a$ =  (1,2)(3,6,7,4,5,8) \quad  $b$ =  (1,7,2,6,4,5)(3,8)  \\
                                                                    \\
    and presentation:                                                 \\
                                                                      \\
  $    a^6=(a*b^{-1})^3=b^6=a^2*b^{-1}*a*b*a^{-2}*b=(a^2*b)^2*a^{-1}*b^{-2}=1  $ \\
                                                                   \\ \hline
                                                                       \\
 	There is one additional case without a normal sylow subgroup,     \\
 	namely $S_4$.                                                        \\
                                                                       \\
 			  Aut($S_4$) = $S_4$.                                     \\
                                                                       \\
 	The automorphism groups of the groups of order $8p$  with a normal   \\
 	sylow $p$-subgroup are obtained by forming the direct product of    \\
 	$Hol(C_p)$ with the entry in the above table.                        \\
 	                                                                  \\ \hline
 	Relations used for the groups of order 8:                         \\
                                                                       \\
 	$C_8:  a^8=1 \qquad   C_4 \times  C_2: a^4=b^2=(a,b)=1 \qquad D_4: a^4=b^2=ab*a=1  $   \\
                                                                       \\
 	  $ C_2 \times  C_2 \times  C_2: a^2=b^2=c^2=(a,b)=(a,c)=(b,c)=1 $     \\
                                                                       \\
 			$Q_2:  a^4=b^4=a^2*b^2=ab*a=1$              \\ \hline
\end{tabular}

\begin{tabular}{|c|l|l|c|} \hline
\multicolumn{4}{|c|}{Table 2            }\\ \hline
\multicolumn{4}{|c|}{Groups of order $8*p^2$ and their automorphism groups;}\\
\multicolumn{4}{|c|}{ cases with a normal sylow $p$-subgroup  }\\ \hline
  	 & \multicolumn{3}{|c|}{ image of 2-group} \\ \cline{2-4}
 &&& \\
2-group&  $C_2\times C_2$ & $C_4$  & order 8 cases\\ 
  & & & \\ \hline
 	     &                    &                  &	                      \\ 

   $C_8$     &                    & [a]\quad $C_2 \times  [144]* $ &  [a]   [144]*            \\
  	     &                  &                  &  for $p \equiv  1$ mod(8)        \\
  	     &                    &                  &     see table 6          \\ \hline
  $C_4 \times  C_2$ & [a,b] $G = Q_p \times  D_p$ & [a]\quad$ G = C_2 \times $     & For $p \equiv  1$ mod(4)         \\
  	     &                    &  $(C_p \times  C_p) @ C_4$ &                          \\
  	     &  $C_2 \times  C_2 \times$ &                  &  [a,b]  $G = Hol(C_p) \times  D_p$\\
  	     & $Hol(C_p) \times  Hol(C_p)$  & $ C_2 \times  [144] * $   & $  Hol(C_p) \times  Hol(C_p)$      \\ \cline{2-4}
  	     & [ab,b]             &                  & For $p \equiv  1$ mod(4)         \\
  	     &                    &                 &                          \\
  	     &$(Hol(C_p) \times  C_2) wr C_2$&                  &  [a,ab] $ Hol(C_p) wr C_2$   \\ \hline
 $ C_2 \times  C_2 $&  [a,b]             &                  &  	                  \\
   $\times  C_2$  &   $G = D_p \times  D_p \times  C_2$&                 &  		           \\
  	     &                    &                  &  	                  \\
  	     & $(Hol(C_p) \times  C_2)wr C_2$&                  &                          \\ \hline
  	     &	 [a,b]        &                  &   for $p = 3$, $G = S_3 wr C_2$ \\
          &                    &                 &     and Aut($G$) is a      \\
          &   $C_2 \times  C_2 \times $      &                  &   is a complete group    \\
     $D_4$   & $ Hol(C_p) \times  Hol(C_p)$ &                  &   of order 144 = [144] * \\ \cline{2-4}
          &    [ab,a]          &                  &  	                  \\
  	     &	              &                  &  	                  \\
  	     &$(Hol(C_p) \times  C_2) wr C_2$&                  &  	                  \\ \hline
 	     &    [ab,a]          &                  &  $p=3$  $Hol(C_3 \times  C_3)$       \\
   $Q_2$     &                    &                 &   	                  \\
 	     &$(Hol(C_p) \times  C_2) wr C_2$&                  &  $p > 3$ complete group    \\
  	     &                    &                  &       not $Hol(C_p \times  C_p)$ * \\ \hline
  totals  &         6          &        2         &  	    3             \\
 	     &                    &                  &                          \\ \hline
\multicolumn{4}{|c|}{ *  See Notes for Table 2 contained in Appendix 1:} \\                                                                                                 
 \multicolumn{4}{|c|}{ for the $D_4$ case see Table D, }    \\
 \multicolumn{4}{|c|}{for the $Q_2$ case see Table Q. } \\ \hline                                 
\end{tabular}

\begin{tabular}{|c|c|c|c|} \hline
\multicolumn{4}{|c|}{ Table 3}\\
\multicolumn{4}{|c|} { } \\
\multicolumn{4}{|c|}{Groups of order $8*p^2$ without a normal $p$-subgroup of odd order}\\ 
\multicolumn{4}{|c|} { } \\ \hline
\multicolumn{4}{|c|} {Groups of order $8*p^2$ with a normal } \\
\multicolumn{4}{|c|} {sylow 2-subgroup  } \\ \hline 
   2-group       &    $p$-group  &       group          &      Aut($G$)       \\ \hline
                 &             &                      &                   \\
 $C_2 \times  C_2 \times  C_2 $ & $  C_3 \times  C_3$   & $  A_4 \times  C_3 \times  C_2$       &    $S_4 \times  S_3$        \\
                 &             &                      &                   \\
              &             &                      &                   \\
      &    $ C_9$      &  $(C_2 \times  C_2) @ C_9 \times  C_2$ & $   S_4 \times  C_3$  \\
                 &             &                      &                   \\
                 &             &                      &                   \\
              &   $C_7 \times  C_7$& $(C_2 \times  C_2 \times  C_2) @ C_7 $ & [168]$ \times  Hol(C_7)$ **\\
                 &             &      $  \times  C_7$          &                   \\
                 &             &                      &                   \\
    & $ C_{49}$  & $(C_2 \times  C_2 \times  C_2)) @ C_{49}$&  [168]$ \times  C_6$     **\\
                 &             &                      &                   \\ \hline
                 &             &                      &                   \\
     $Q_2$     & $ C_3 \times  C_3$  & $SL(2,3) \times  C_3$  & $ S_4 \times  S_3$      \\
                 &             &                      &                   \\
                 &  $C_9$     & $Q_2 @ C_9$              &   $ S_4 \times  C_3$        \\
                 &             &                      &                   \\ \hline                                                                          
\multicolumn{4}{|c|} {** [168] is the complete group of this order; see Table 1 for details.}\\ \hline 
\multicolumn{4}{|c|} {Non-normal sylow subgroup types:  $p = 3$.} \\
\multicolumn{4}{|c|} { } \\
\multicolumn{4}{|c|} {$S_4 \times  C_3 \longrightarrow  \text{Aut}(G) = S_4 \times  C_2$  } \\
\multicolumn{4}{|c|} {$A_4 \times  S_3\longrightarrow  \text{Aut}(G) = S_4 \times  S_3 $ } \\
\multicolumn{4}{|c|} {$(C_2 \times  C_2) @ D_9 \rightarrow $ Aut($G$) = complete group of order 216    } \\
\multicolumn{4}{|c|} {$(A_4 \times  C_3 ) @ C_2$	 [even permutations in $S_3 \times  S_4]$.    } \\
\multicolumn{4}{|c|} { Aut($G$) = complete group of order 432     } \\
\multicolumn{4}{|c|} {with 20 conjugacy classes.       } \\
\multicolumn{4}{|c|} { } \\ \hline
\multicolumn{4}{|l|} {A presentation for the complete group of order 216 is:    } \\
\multicolumn{4}{|c|} { } \\
\multicolumn{4}{|c|} {$a^2=b^2=c^3=(a,c)=(a*d)^2=(b,c)=c*d^2*(c^{-1})*d$     }\\
\multicolumn{4}{|c|} {$= a*b*a*d*b*(d^{-1})=(a*b)^2*(d^{-1})*b*d=1;   $          } \\
\multicolumn{4}{|c|} { } \\
\multicolumn{4}{|l|} {A presentation for the complete group of order 432 is:      } \\
\multicolumn{4}{|c|} { } \\
\multicolumn{4}{|c|} {$a^4=b^2=c^3=d^2=(a,c)=(a,d)=(a,d)=(b,d)=(c*d)^2=(a*b)^3 $    } \\
\multicolumn{4}{|c|} {$=(b*c)^2*(b*(c^{-1}))^2=a*b*a*c*b*c*a*b*c=1; $           } \\
\multicolumn{4}{|c|} { } \\ \hline
\end{tabular}

\begin{tabular}{|l|c|c|c|} \hline
\multicolumn{4}{|c|}{} \\
\multicolumn{4}{|c|}{ Table 4} \\
\multicolumn{4}{|c|}{} \\
\multicolumn{4}{|c|}{$C_8$ extensions from Yuanda's work [16]}\\
\multicolumn{4}{|c|}{} \\ \hline
        prime	  &     $C_2$ image	 &     $ C_4$ image &     $ C_8 $ image \\ \hline
                    &                      &                   &     \\
         $p \equiv 1$ mod(8) &      (2,6)     &      (3,5,8,9)  &   (4,7,10,11,12, \\
    			&			 &	              &     13,14,15)  \\
                  &                 &                  &                  \\
         $p \equiv 3$ mod(8)&       (2,6)    &         (5)      &  (16)           \\
                    &               &                  &                 \\
         $p \equiv 5$ mod(8)&	   (2,6)     &    (3,5,8,9)     &(7)            \\
                    &               &                  &               \\
         $p \equiv 7$ mod(8)&	   (2,6)     &     (5)          &(17)           \\
                     &                &                 &     \\ \hline
\end{tabular}
\pagebreak

\begin{tabular}{|l|c|c|l|} \hline
\multicolumn{4}{|c|}{} \\
\multicolumn{4}{|c|}{Table 5} \\
\multicolumn{4}{|c|}{Explicit representations for $(C_p \times  C_p) @ C_8$} \\
\multicolumn{4}{|c|}{for some small primes} \\ 
\multicolumn{4}{|c|}{} \\ \hline
        prime	&	$C_2$ image	   &  $C_4$ image	 &\quad   $C_8$ image  \\ \hline
                  &                    &               &                \\
                  &                    &               &                \\
       $p=3$       &    (1,0;0,-1)      &(0,-1;1,0)     & 16)=(0,1;1,1)  \\
  	            &    (1,0;0,1)       &               &                \\ \hline
       $p=5$       &    (1,0;0,-1)      &(2,0;0,-1)      & 7)=(0,1;2,0)     \\
                 &    (1,0;0,1)       &(2,0;0,1)       &                 \\
        		&	               &(2,0;0,2)       &                  \\
        		&                    &(2,0;0,3)       &                \\ \hline
       $p=7$       &    (1,0;0,-1)      &(0,-1;1,0)     & 17)=(0,1;-1,3)   \\
        	      &    (1,0;0,1)       &                &                \\ \hline
       $p \equiv 1$ mod(8)&    (1,0;0,-1)      &(4,0;0,-1)      & 4)=(-2,0;0,-2) \\
       ($p=17$).     &   (1,0;0,1)       &(4,0;0,1)       & 7)=(0,-1;4,0)   \\
                   &                  &(4,0,0,4)      & 10)=(-2,0;0,-4)   \\
        		 &	               &(4,0;0,-4)     & 11)=(-2,0;0,4)    \\
                  &                   &               & 12)=(-2,0;0,8)    \\
                  &                   &               & 13)=(-2,0;0,-8)   \\
                  &                   &               & 14)=(-2,0;0,1)    \\
                  &                   &               & 15)=(-2,0;0,-1) * \\
                  &                   &                &            \\ \hline
 \multicolumn{4}{|c|}{  }\\                                                           
 \multicolumn{4}{|c|}{* Zhang Yuanda says number 15 has a center of order 1.}      \\
 \multicolumn{4}{|c|}{  This is not the case. This group has the form }          \\
 \multicolumn{4}{|c|}{ $ ( C_{17} @ C_{16} ) \times C_{17} $.  }      \\ \hline                                                  
\end{tabular}

\begin{tabular}{|l|c|c|c|} \hline
\multicolumn{4}{|c|}{Table 6} \\
\multicolumn{4}{|c|}{Automorphism groups for the groups given in Table 5 }\\ \hline
      prime	& $  C_2$ image	&       $ C_4$ image   &    $  C_8$ image  \\ \hline
            &                 &                      &                   \\
    $p=3$ & $ Hol(C_3) \times  C_2$& $ H(3^2) \times C_2 $  & $ 16) = H(3^2) $  \\
       &$\times C_2 \times C_2$  & & \\ \cline{2-2}
  	  & $Hol(C_3 \times C_3) $   & & \\ 
        & $\times C_2 \times C_2$ & & \\ \hline
    $p=5$ &  $Hol(C_5) \times C_4 $ & $ Hol(C_5) \times C_4 \times C_2$ &  7) = $H(5^2)$     \\
         &   $ \times C_2\times C_2$ &$Hol(C_5) \times Hol(C_5) \times C_2$    & \\ \cline{2-2}           
        &  $Hol(C_5 \times C_5) $ &   $ Hol(C_5 \times C_5) \times C_2$ &\\
        & $\times C_2\times C_2$ &$(Hol(C_5) wr C_2) \times C_2$        & \\   \hline
    $p=7$  & $Hol(C_7) \times C_6  $ &  $  H(7^2) \times C_2 $&    17) = $H(7^2)$      \\
        & $\times C_2\times C_2$& & \\ \cline{2-2}
         & $Hol(C_7 \times C_7)  $ &   &  \\
          & $ \times C_2\times C_2$ & & \\ \hline
   $p \equiv 1$ mod(8)& \multicolumn{2}{|c|}{}       &    4)=$Hol(C_p \times C_p)$     \\
   (e.g.,$p=17$). & \multicolumn{2}{|c|}{See comments below for the}& 7)=$Hol(C_p) wr C_2$    \\
                & \multicolumn{2}{|c|}{$C_2$ and $C_4$ image cases.}&       10)=$Hol(C_p) \times Hol(C_p)$ \\
                & \multicolumn{2}{|c|}{}    &       11)=$Hol(C_p) \times Hol(C_p)$ \\ 
                & \multicolumn{2}{|c|}{}    &       12)=$Hol(C_p) wr C_2$    \\
                & \multicolumn{2}{|c|}{}    &          13)=$Hol(C_p) wr C_2 $  \\
                & \multicolumn{2}{|c|}{}    &          14)=$Hol(C_p) \times Hol(C_p)$ \\
                & \multicolumn{2}{|c|}{}    &          15)=$Hol(C_p) \times C_{(p-1)}$ \\ \hline
\multicolumn{4}{|c|}{$H$($p^n$) is the group of all mappings $x \mapsto  a*x^t + b,$}\\
\multicolumn{4}{|c|}{where $a \not= 0$ and $b$ are elements of the Galois group of}\\   
\multicolumn{4}{|c|}{$F$ = GF($p^n$) and $t$ is a field automorphism of $F$.}\\                     
\multicolumn{4}{|c|}{The order of $H$($p^n$) \text{ is} $p^n * ( p^n  -  1) * n.$} \\ \hline 
\multicolumn{4}{|c|}{For the $p \equiv 3$ mod(8) and $p \equiv 7$ mod(8) cases, the entries in }\\    
\multicolumn{4}{|c|}{Table 6 have the following generalizations:     }\\
\multicolumn{4}{|c|}{In the $C_2$ image cases just make the following substitutions:}\\
\multicolumn{4}{|c|}{ $Hol(C_3) \times C_2$ \text{by} $Hol(C_p) \times C_{(p-1)}$ and} \\                   
\multicolumn{4}{|c|}{ $Hol(C_3 \times C_3)$ \text{by} $Hol(C_p \times Cp). $}\\
\multicolumn{4}{|c|}{In the $C_4$ image cases just make the following substitutions:}\\  
\multicolumn{4}{|c|}{[144] (Table 2)  by Aut[ ($C_p \times C_p) @ C_4] = H(p^2)$. }\\
\multicolumn{4}{|c|}{The analogue, or recurrence, for the $p \equiv 1$ mod(4)  }\\
\multicolumn{4}{|c|}{ case is  $Hol(C_p) wr C_2. $} \\              
\multicolumn{4}{|c|}{In the $C_8$ image cases just make the following substitutions:}\\  
\multicolumn{4}{|c|}{ $H(3^2)$ [ $H(7^2)$ respectively ] by $ H(p^2)$. }\\ \hline            
\multicolumn{4}{|c|}{For $p \equiv 1, 5$ mod(8) make the following changes in the $p$ = 5 cases: }\\
\multicolumn{4}{|c|}{In the $C_2$ image cases just make the following substitutions:}\\  
\multicolumn{4}{|c|}{ $C_5$ by $C_p$ and $C_4$ by $C_{(p-1)} . $}\\
\multicolumn{4}{|c|}{In the $C_4$ image cases just make the following substitutions:} \\ 
\multicolumn{4}{|c|}{$C_5$ by $C_p$ and $C_4$ by $C_{(p-1)}$.}\\                             
\multicolumn{4}{|c|}{In the $C_8$ image cases just make the following substitutions:}\\  
\multicolumn{4}{|c|}{ $H(5^2)$  by $H(p^2)$ [$p \equiv 5$ mod(8) here]. }\\ \hline
\end{tabular}

\begin{thebibliography}{99}
\bibitem{3} In addition to these large general purpose systems there are a
      number of smaller more teaching-oriented programs for 
	teaching abstract algebra. The groups discussed here are not 
	of very high order and may be able to be run on many of 
	these smaller educationally oriented group theory systems.\\
\bibitem{5} The book by Fuchs on \textquotedblleft Abelian Groups" is one of the few exceptions.
        This book does have a few pages devoted to the discussion of 
	the automorphism groups of abelian groups. This book also 
	contains references to the articles of K. Shoda on the automorphism
	groups of abelian groups. See also the brief 
	discussion on automorphism groups of abelian groups in 
	Carmichael's textbook on finite groups [6].\\
\bibitem{12}  The wreath product of the group $T$ with $C_p$  ($T wr C_p$) can be 
	described as follows. Write the generators of the group $T$
	in terms of permutations of some degree, say $n$. Then the
	group $C_p$ is expressed in the form:
		$C_p$: (1,$n+1,2n+1,...,(p-1)*n + 1$)
		    (2,$n+2,2n+2,...,(p-1)*n + 2$)
		     $ \cdots$
		    ($n, 2n, 3n, ...,p*n$).\\

      That is, we have $n$ cycles of length $p$ in this representation.   
      In the particular cases below of $(Hol(C_p) \times  C_2) wr C_2$ we have 
      for $p = 3$ the following:
			$a$ = (1,2,3), \quad  $b$ = (2,3), \quad $c$ = (4,5), \quad 
			   $d$ = (1,6)(2,7)(3,8)(4,9)(5,10).

      Here the permutations $a$, $b$ and $c$ generate the group $Hol(C_3) \times  C_2$.\\
\bibitem{15} The way to determine the exact number of cases here is to
	look at the normal subgroup lattices for the groups of 
	order 64 given in the Hall-Senior book [11] and actually
      count the number of different $C_2$ quotient groups listed 
	there. There are 1120 such cases; 700 of these are 
	associated with a characteristic subgroup. As one can readily
	see this is a rather tedious undertaking.\\
\bibitem{17} To do this for an arbitrary group (i.e., to give a finite
	list of $n$ quantities, no two sets of which are the same for 
	nonisomorphic groups) appears to be a very hard problem in 
	group theory. \\
\bibitem{20} Relations of the form $x^t \equiv 1$ mod($p$) occur in several places
	in the text and the appendix. In each case, $x$ (or $q$ or $z$) is 
	a $t$-th root of unity (here $t$ = 3, 4, ..., ($p-1$) ).  \\
\end{thebibliography}
\end{document}